\documentclass[12pt]{amsart}

\usepackage{amstext}
\usepackage{amsthm}
\usepackage{amsrefs}
\usepackage{amsmath}
\usepackage{amssymb}
\usepackage{latexsym}
\usepackage{amsfonts}
\usepackage{graphicx}
\usepackage{texdraw}
\usepackage{graphpap}

\input txdtools

\begin{document}

\newcommand{\ci}[1]{_{ {}_{\scriptstyle #1}}}

\newcommand{\norm}[1]{\ensuremath{\|#1\|}}
\newcommand{\abs}[1]{\ensuremath{\vert#1\vert}}
\newcommand{\p}{\ensuremath{\partial}}
\newcommand{\pr}{\mathcal{P}}

\newcommand{\pbar}{\ensuremath{\bar{\partial}}}
\newcommand{\db}{\overline\partial}
\newcommand{\D}{\mathbb{D}}
\newcommand{\B}{\mathbb{B}}
\newcommand{\Sp}{\mathbb{S}}
\newcommand{\T}{\mathbb{T}}
\newcommand{\R}{\mathbb{R}}
\newcommand{\Z}{\mathbb{Z}}
\newcommand{\C}{\mathbb{C}}
\newcommand{\N}{\mathbb{N}}
\newcommand{\scrH}{\mathcal{H}}
\newcommand{\scrL}{\mathcal{L}}
\newcommand{\td}{\widetilde\Delta}

\newcommand{\La}{\langle }
\newcommand{\Ra}{\rangle }
\newcommand{\rk}{\operatorname{rk}}
\newcommand{\card}{\operatorname{card}}
\newcommand{\ran}{\operatorname{Ran}}
\newcommand{\im}{\operatorname{Im}}
\newcommand{\re}{\operatorname{Re}}
\newcommand{\tr}{\operatorname{tr}}
\newcommand{\vf}{\varphi}
\newcommand{\f}[2]{\ensuremath{\frac{#1}{#2}}}


\newcommand{\entrylabel}[1]{\mbox{#1}\hfill}

\newenvironment{entry}
{\begin{list}{X}%
  {\renewcommand{\makelabel}{\entrylabel}%
      \setlength{\labelwidth}{55pt}%
      \setlength{\leftmargin}{\labelwidth}
      \addtolength{\leftmargin}{\labelsep}%
   }%
}%
{\end{list}}



\numberwithin{equation}{section}

\newtheorem{thm}{Theorem}[section]
\newtheorem{lm}[thm]{Lemma}
\newtheorem{cor}[thm]{Corollary}
\newtheorem{conj}[thm]{Conjecture}
\newtheorem{prob}[thm]{Problem}
\newtheorem{prop}[thm]{Proposition}
\newtheorem*{prop*}{Proposition}

\theoremstyle{remark}
\newtheorem{rem}[thm]{Remark}
\newtheorem*{rem*}{Remark}

\title{A Note about Stabilization in $A_\R(\D)$}
\author{Brett D. Wick}
\address{Brett D. Wick, Department of Mathematics\\ Vanderbilt University\\ 1326 Stevenson Center\\ Nashville, TN USA 37240-0001}
\email{brett.d.wick@vanderbilt.edu}
\thanks{Research supported in part by a National Science Foundation RTG Grant to Vanderbilt University}



\begin{abstract} 
It is shown that for $A_\R(\D)$ functions $f_1$ and $f_2$ with 
$$
\inf_{z\in\overline{\D}}(\abs{f_1(z)}+\abs{f_2(z)})\geq\delta>0
$$
and $f_1$ being positive on real zeros of $f_2$ then there exists $A_\R(\D)$ functions $g_2$ and $g_1,g_1^{-1}$ with and 
$$
g_1f_1+g_2f_2=1\quad\forall z\in\overline{\D}.
$$


This result is connected to the computation of the stable rank of the algebra $A_\R(\D)$ and to Control Theory.
\end{abstract}

\maketitle

\section{Introduction and main result}
\label{Intro.}
The stable rank of a ring (also called the Bass stable rank) was introduced by H. Bass in \cite{Bass} to assist in computations of algebraic K-Theory.  We recall the definition of the stable rank.

Let $A$ be an algebra with a unit $e$.  An n-tuple $a\in A^n$ is called \textit{unimodular} if there exists an n-tuple $b\in A^n$ such that $\sum_{j=1}^{n}a_jb_j=e$.  A n-tuple $a$ is called \textit{stable} or \textit{reducible} if there exists an (n-1)-tuple $x$ such that the (n-1)-tuple $(a_1+x_1a_n,\ldots, a_{n-1}+x_{n-1}a_{n})$ is unimodular.  The \textit{stable rank} (also called \textit{bsr}(A) in the literature) of the algebra $A$ is the least integer n such that every unimodular (n+1)-tuple is reducible.  

The stable rank is a purely algebraic concept, but can be combined with analysis when studying commutative Banach algebras of functions.  In this context the stable rank is related to the zero sets of ideals of the Banach algebra and the spectrum of the Banach algebra.  The stable rank for different algebras of analytic functions have been considered by many authors.  The computation of the stable rank of the disc algebra $A(\D)$ was done by Jones, Marshall and Wolff, \cite{JMW}.  Their main result was the following theorem
\begin{thm}[P. Jones, D. Marshall and T. Wolff, \cite{JMW}]
\label{stabledisc}
Let $f_1,f_2\in A^\infty(\D)$ be such that $\inf_{z\in\overline{\D}}(\abs{f_1(z)}+\abs{f_2(z)})=\delta>0$.  Then there exists $g_1, g_2, g_1^{-1}\in A^\infty(\D)$ with   
$$
f_1(z)g_1(z)+f_2(z)g_2(z)=1\quad\forall z\in\overline{\D}.
$$
\end{thm}
It is immediate that this result implies that stable rank of $A(\D)$ is one.  The computation was done for sub-algebras of the disk algebra $A(\D)$ by Corach--Su\'arez, \cite{CS}, and Rupp \cite{R}.  This paper concerns results related to the stable rank of the algebra $A_\R(\D)$, a certain sub-algebra of the disc algebra.

It is possible to phrase the result \cite{JMW} in the language of Control Theory.  In this language, the result can be viewed as saying that it is possible to stabilize (\textit{in the sense given above}) a linear system (\textit{the Corona data, viewed as a rational function}) via a stable (\textit{analytic}) controller.  But, in applications of Control Theory, the linear systems and transfer functions  have real coefficients, so in this context Theorem \ref{stabledisc} is physically meaningless.  From the point of view of Control Theory, it is important to know if results like Theorem \ref{stabledisc} hold, but for a more physically meaningful algebra, and is the main motivation for this paper.  This paper is interested in questions related to the stable rank of a natural sub-algebra of $A(\D)$, the real Banach algebra $A_\R(\D)$, in particular, does some variant of Theorem \ref{stabledisc} hold for this algebra. 

First, recall that $A_\R(\D)$ is the subset of $A(\D)$ with additional property that the Fourier coefficients of an element of $A_\R(\D)$ must be real.  This property can be captured by the following symmetry condition
$$
f(z)=\overline{f(\overline{z})}\quad\forall z\in\overline{\D}.
$$
This condition is implying that the functions in $A_\R$ possess a symmetry that is not present  for general $A$ functions.

The Corona result is inherited by the algebra $A_\R(\D)$.  More precisely, it is an immediate application of the usual Corona Theorem, see \cite{Rudin}, and the symmetry properties of $A_\R(\D)$ to show that an 
n-tuple $(f_1,\ldots,f_n)\in (A_\R(\D))^n$ is unimodular if and only if it satisfies the Corona condition,
$$
\inf_{z\in\overline{\D}}\left(\abs{f_1(z)}+\cdots+\abs{f_n(z)}\right)=\delta>0.
$$
Indeed, one direction is immediate, and in the other direction, if we know that
$$
\inf_{z\in\overline{\D}}\left(\abs{f_1(z)}+\cdots+\abs{f_n(z)}\right)=\delta>0,
$$
then we can find a solution $(g_1,\ldots,g_n)\in (A(\D))^n$.  We then symmetrize the $g_j$ via the operation
$$
\tilde{g}_j(z):=\f{g_j(z)+\overline{g_j(\overline{z})}}{2}.
$$
The $\tilde{g}_j\in A_\R(\D)$ and will then be the $A_\R(\D)$ Corona solution we are seeking.

This leads to the main question considered in this paper.  Is Theorem \ref{stabledisc} true for the algebra $A_\R(\D)$? Namely, given Corona data $f_1$ and $f_2$ in $A_\R(\D)$ is there a solution $g_1$ and $g_2$ to the Corona problem with $g_1$ invertible in $A_\R(\D)$?  

It is very easy to see that there is an additional necessary condition.  Suppose that Theorem \ref{stabledisc} were true for $A_\R$ functions, then we shall see that the real zeros of $f_1$ and $f_2$ must intertwine correctly.  Indeed, let $\lambda_1$ and $\lambda_2$ be real zeros of $f_2$.  Then we have
\begin{eqnarray*}
f_1(\lambda_1)g_1(\lambda_1) & = & 1\\
f_1(\lambda_2)g_1(\lambda_2) & = & 1.
\end{eqnarray*}
Now $f_1(\lambda_1)$ and $f_1(\lambda_2)$ must have the same sign at these zeros.  If this were not true then without loss of generality suppose that $f_1(\lambda_1)>0>f_1(\lambda_2)$.  Then $g_1(\lambda_1)>0>g_1(\lambda_2)$.  By continuity there will exists a point $\lambda_{12}$, between $\lambda_1$ and $\lambda_2$, with $g_1(\lambda_{12})=0$.  But this contradicts the fact that $g_1^{-1}\in A_\R(\D)$.  So $f_1$ must have the same sign at real zeros of $f_2$.  We will say that $f_1$ is POZ of $f_2$, if $f_1$ has the same sign at all real zeros of $f_2$.  Here POZ stands for \textit{positive on real zeros}.

This is also an intertwining condition of the zeros of $f_1$ and $f_2$.  More precisely, the functions $f_1$ and $f_2$ satisfy condition POZ if and only if between every real zero of $f_2$ there must be an even number of real zeros of $f_1$.  This is called the \textit{parity interlacing property} found in Control Theory for the stabilization of a linear system.

The main result of this paper is the following theorem.

\begin{thm}
\label{stablerealdisc}
Suppose that $f_1,f_2\in A_\R(\D)$, $\norm{f_1}_\infty,\norm{f_2}_{\infty}\leq 1$, $f_1$ is POZ of $f_2$ and
$$
\inf_{z\in\overline{\D}}\left(\abs{f_1(z)}+\abs{f_2(z)}\right)=\delta>0.
$$
Then there exists $g_1, g_1^{-1}, g_2\in A_\R(\D)$ with 
$$
f_1(z)g_1(z)+f_2(z)g_2(z)=1\quad\forall z\in\overline{\D}.
$$
\end{thm}



Theorem \ref{stablerealdisc} also immediately implies the that the stable rank of $A_\R(\D)$ is at least 2.

We remark that these results transfer immediately to analogous statements $A_\R(\C_+)$ via the standard conformal mapping between $\C_+$ and $\D$.

Throughout the paper, the adjectives \textit{real symmetric} are used to indicate that the function in question satisfies the symmetry condition 
$$
f(z)=\overline{f(\overline{z})}\quad\forall z\in\overline{\D}.
$$

\section{Proof of Main Result}
\label{Main}

The proof of the main result is inspired by the proof used in Jones, Marshall and Wolff, \cite{JMW}.  We use the symmetry of $A_\R(\D)$ functions and the condition that $f_1$ is POZ of $f_2$ to modify the proof in \cite{JMW}.  

\subsection{The Key Proposition}
The proof can be deduced from the following proposition.

\begin{prop}
\label{LipFunction}
Suppose that $f_1,f_2\in A_\R(\D)$, $\norm{f_1}_\infty,\norm{f_2}_{\infty}\leq 1$, $f_1$ is POZ of $f_2$ and
$$
\inf_{z\in\overline{\D}}\left(\abs{f_1(z)}+\abs{f_2(z)}\right)=\delta>0.
$$
Then there exists a continuous Lipschitz function $\Phi:\overline{\D}\to\C$ and $\delta$ with $0<\delta$ such that
\begin{itemize}
\item[(1.)] $\Phi(z)=f_1(z)$ if $\abs{f_1(z)}<\delta$;
\item[(2.)] $\Phi(z)=1$ if $\abs{f_2(z)}<\delta$;
\item[(3.)] $\abs{\Phi(z)}\geq\delta$ if $\abs{f_1(z)}\geq\delta$;
\item[(4.)] $\f{\p \Phi}{\p\overline{z}}$ is a bounded function on $\D$;
\item[(5.)] $\Phi(z)=\overline{\Phi(\overline{z})}$ for all $z\in\overline{\D}$.
\end{itemize}
\end{prop}

\begin{proof}
To construct $\Phi$, choose closed sets $E_1$ and $E_2$ such that $E_j$ has finitely many components and 
$$
\{z: \abs{f_j(z)}<\delta\}\subseteq E_j\subseteq\{z:\abs{f_j}<2\delta\}
$$
where $\delta$ is chosen so that 
$$
\{z:\abs{f_1(z)}<2\delta\}\cap\{z:\abs{f_2(z)}<2\delta\}=\emptyset.
$$
Using the symmetry property of $A_\R(\D)$ functions, it is possible to select the sets $E_j$ so that the are symmetric with respect to the real axis, namely, $E_j=\overline {E_j}$ (with bar denoting conjugation).  By the maximum principle, $E_2$ can not separate any point of $E_1$ from the boundary of the disc.  We can then extend the components of $E_1$ to the boundary to obtain a closed set $S\subset\D$ with the following properties
\begin{itemize}
\item[(a)] $E_1\subset S$;
\item[(b)] $S\cap E_2=\emptyset$;
\item[(c)] $S$ has finitely many components, each of which intersects the boundary of the disc;
\item[(d)] $S$ is symmetric with respect to the real axis, i.e., $S=\overline{S}$.
\end{itemize}
The set $\overline{\D}\setminus S$ consists of simply connected components.  Since $f_1$ is POZ of $f_2$ there will exist a bounded real symmetric branch of $\log f_1$ on $\overline{\D}\setminus S$, see \cite{Vidy}.  Since $E_2$ and $S$ are compact and symmetric, and since $E_2\cap S=\emptyset$, there will exist a function $q(z)$ with the properties,
\begin{itemize}
\item[(a)] $q(z)=1$ on a neighborhood of $S$;
\item[(b)] $q(z)=0$ on a neighborhood of $E_2$;
\item[(c)] $q(z)=\overline{q(\overline{z})}$ for all $z\in\overline{\D}$.
\end{itemize}

With this function we then define $\Phi(z):=\exp(q(z)\log f_1(z))$.  It is then straightforward to demonstrate that this function satisfies the conclusions of Proposition \ref{LipFunction}.
\end{proof}

\subsection{Proof of Theorem \ref{stablerealdisc}}
Set $g_1(z):=\f{\Phi(z)}{f_1(z)}\exp(u(z) f_2(z))$ for some function $u$ to be defined.  If $u$ is real symmetric,  analytic on $\D$, continuous on $\overline{\D}$, and
$$
\f{\p u }{\p\overline{z}}= \f{1}{\Phi f_2}\f{\p \Phi}{\p\overline{z}}
$$
then $g_1\in A_\R(\D)$.  Here we use Proposition \ref{LipFunction}.

Let 
$$
k:=\f{1}{\Phi f_2}\f{\p \Phi}{\p\overline{z}}.
$$
Using Proposition \ref{LipFunction} we have that $k$ is a bounded real symmetric function on $\D$.  Finally, set
$$
u(z):=\f{1}{\pi}\int_\D\f{k(\xi)}{\xi-z}dA(\xi)
$$
Then $u$ is real symmetric (by direct verification), continuous on $\overline{\D}$ (since it is the convolution of a bounded function with a integrable function), and satisfies $\f{\p u}{\p\overline{z}}=k$, thus $g_1\in A_\R(\D)$.  Moreover, we also have that $g_1$ is bounded away from zero, so $g_1^{-1}\in A_\R(\D)$.

Finally, let $g_2$ be defined by the condition $f_1g_1+f_2g_2=1$.  Clearly, $g_2$ is continuous when $f_2\neq 0$.  But, if $0<\abs{f_2}<\delta$ then 
$$
g_2=\f{1-f_1g_1}{f_2}=\f{1-\exp(uf_2)}{f_2}
$$
by Proposition \ref{LipFunction}.  Finally, note that $g_2\to -u$ as $f_2\to 0$.  This finishes the proof of Theorem \ref{stablerealdisc}.

\section{Concluding Remarks}
\label{remarks}

The argument above shows that the stable rank of $A_\R(\D)$ is at least two.  A natural conjecture is

\begin{conj}
\label{HardConj}
The stable rank of $A_\R(\D)$ is two.
\end{conj}

To solve this conjecture one needs to answer the following question:

\begin{prob}
\label{Hard}
Let $f_1, f_2, f_3\in A_{\R}(\D)$ be such that $\inf_{z\in\overline{\D}}(\abs{f_1(z)}+\abs{f_2(z)}+\abs{f_3(z)})=\delta>0$.  Then does there exist $h_1, h_2, g_1, g_2\in A_{\R}(\D)$ such that
$$
f_1(z)g_1(z)+f_2(z)g_2(z)+(h_1(z)g_1(z)+h_2(z)g_2(z))f_3(z)=1\quad\forall z\in\overline{\D}?
$$
\end{prob}

The stable rank being two is simply saying that it is possible to choose a Corona solution with the property that one of the solution elements is in the ideal generated by the others.  This property is true in the algebra $A(\D)$, but it is not clear how to force this property onto the sub-algebra $A_\R(\D)$.  

Additionally, we remark because of the additional property of continuity at the boundary of $A_\R(\D)$ functions, many of the analytic difficulties in the proof of Theorem \ref{stablerealdisc} .  Similar results can be demonstrated for the algebra $H^\infty_\R(\D)$, the sub-algebra of $H^\infty(\D)$ with real Fourier coefficients.  The method of proof requires different techniques, and can be found in \cite{Wick}.


\section*{References}

\begin{biblist}
\bib{Bass}{article}{AUTHOR = {Bass, Hyman},
     TITLE = {{$K$}-theory and stable algebra},
   JOURNAL = {Inst. Hautes \'Etudes Sci. Publ. Math.},
    NUMBER = {22},
      YEAR = {1964},
     PAGES = {5--60}
}


\bib{CS}{article}{
   author={Corach, Gustavo},
   author={Su{\'a}rez, Fernando Daniel},
   title={Stable rank in holomorphic function algebras},
   journal={Illinois J. Math.},
   volume={29},
   date={1985},
   number={4},
   pages={627--639},
}


\bib{JMW}{article}{
   author={Jones, P. W.},
   author={Marshall, D.},
   author={Wolff, T.},
   title={Stable rank of the disc algebra},
   journal={Proc. Amer. Math. Soc.},
   volume={96},
   date={1986},
   number={4},
   pages={603--604},
   }
   
   
\bib{Rudin}{book}{
   author={Rudin, Walter},
   title={Real and complex analysis},
   edition={3},
   publisher={McGraw-Hill Book Co.},
   place={New York},
   date={1987},
   pages={xiv+416},
   isbn={0-07-054234-1},
}
   
\bib{R}{article}{
   author={Rupp, Rudolf},
   title={Stable rank of subalgebras of the disc algebra},
   journal={Proc. Amer. Math. Soc.},
   volume={108},
   date={1990},
   number={1},
   pages={137--142},
}

\bib{Vidy}{book}{
    author={Vidyasagar, M.},
     title={Control System Synthesis: A Factorization Approach},
    series={MIT Press Series in Signal Processing, Optimization, and
            Control, 7},
      publisher={MIT Press},
     place={Cambridge, MA},
      date={1985},
     pages={xiii+436},
      isbn={0-262-22027-X},
 }

\bib{Wick}{article}{
     author={Wick, Brett D.},
   title={Stablization in $H^\infty_\R(\D)$},
   journal={submitted to Math. Ann.}
}
\end{biblist}
\end{document}